\documentclass[a4paper,12pt, preprint, 3p, times]{P-article}
\usepackage{amssymb, amsmath,amsfonts, amsthm, amscd, epsfig, graphics}

\usepackage[english]{babel} 
\usepackage[latin1]{inputenc} 
\usepackage[T1]{fontenc} 

\usepackage{euscript}
\usepackage{wasysym}
\usepackage{ulem}
\usepackage{xcolor}
\usepackage{yfonts, mathrsfs, eufrak}
\usepackage{hyperref}
\usepackage{cleveref}

\usepackage{natbib}

\DeclareMathAlphabet{\mathpzc}{OT1}{pzc}{m}{it}

\newtheorem{thm}{Theorem}[section]
\newtheorem{lem}[thm]{Lemma}

\newtheorem{cor}[thm]{Corollary}

\newtheorem{defn}[thm]{Definition}
\newtheorem{ex}[thm]{Example}

\newtheorem{rem}[thm]{Remark}

\newcommand{\comment}[1]{}

\newcommand{\A}[1]{\mathbb{A}^{#1}}
\newcommand{\KerD}{\text{Ker}(D)}
\newcommand{\Ker}[1]{\text{Ker}(#1)}
\newcommand{\Sym}[2]{\text{Sym}_{#1}(#2)}

\newcommand{\ResRk}[1]{\text{Res-Rk}(#1)}
\newcommand{\ResVarRk}[1]{\text{ResVar-Rk}(#1)}
\newcommand{\Rk}[1]{\text{Rk}(#1)}

\newcommand{\trdeg}[2]{\text{tr.deg}_{#2}(#1)}

\newcommand{\ul}[1] {\underline{#1}}

\newcommand{\Spec}[1]{\text{Spec}(#1)}
\newcommand{\mbbQ}{\mathbb{Q}}

\newcommand{\mbbC}{\mathbb{C}}

\newcommand{\Qt}[1]{Qt(#1)}

\begin{document}
	
	\begin{frontmatter}
		\title{Rank and rigidity of \\
			locally nilpotent derivations of affine fibrations}
		\author{Janaki Raman Babu}
		\address{Department of Mathematics, Indian Institute of Space Science and Technology, \\
			Valiamala P.O., Trivandrum 695 547, India\\
			email: \texttt{raman.janaki93@gmail.com, janakiramanb.16@res.iist.ac.in}}
		
		\author{Prosenjit Das}
		\address{Department of Mathematics, Indian Institute of Space Science and Technology, \\
			Valiamala P.O., Trivandrum 695 547, India\\
			email: \texttt{prosenjit.das@gmail.com, prosenjit.das@iist.ac.in}}
		
		\author{Swapnil A. Lokhande}
		\address{Indian Institute of Information Technology Vadodara,\\
			Block 9, Government Engineering College, \\
			Sector - 28, Gandhinagar 382028, India.\\
			email: \texttt{swapnil@iiitvadodara.ac.in}}

		\begin{abstract}
			In this exposition, we propose a notion of rank and rigidity of locally nilpotent derivations of affine fibrations. We show that the concept is analogous to the perception of rank and rigidity of locally nilpotent derivations of polynomial algebras. Our results characterize locally nilpotent derivations of $\mathbb{A}^3$-fibrations having slice by classifying the fixed point free locally nilpotent derivations in terms of their ranks.\\
			
			\noindent
			\tiny \textbf{Subject Class[2010]}: 13N15, 14R25(primary), 14R20(secondary)\\
			
			\noindent
			\tiny \textbf{Keywords:} Affine fibration, Residual system, Residual variable, Locally nilpotent derivation, Rank, Rigidity.
		\end{abstract}
	\end{frontmatter}
	
	\section{Introduction} \label{Sec_Intro}
	Throughout this article, rings will be commutative with unity. Let $R$ be a ring and $R^{[n]}$ denote the \textit{polynomial ring} in $n$ variables over $R$. Suppose that $A$ is an $R$-algebra. We shall use the notation $A = R^{[n]}$ to mean that $A$ is isomorphic, as an $R$-algebra, to a polynomial ring in $n$ variables over $R$. For a prime ideal $P$ of $R$, $k(P)$ will denote the \textit{residue field} $R_P/PR_P$. If $R$ is a domain, then the notation $\Qt{R}$ will denote the \textit{quotient field} of $R$. $A$ is called an \textit{$\A{n}$-fibration} or \textit{affine $n$-fibration} over $R$, if it is finitely generated and flat over $R$, and $A \otimes_R k(P) = k(P)^{[n]}$ for all $P \in \Spec{R}$. $A$ is said to be a \textit{trivial} $\A{n}$-fibration over $R$ if $A = R^{[n]}$. Let $D: A \longrightarrow A$ be an $R$-derivation. $D$ is called \textit{irreducible} if there does not exist $\alpha \in A\backslash A^*$ such that $D(A) \subseteq \alpha A$. $D$ is defined to be \textit{fixed point free}, if $D(A)A = A$. $D$ is said to have a \textit{slice} $s \in A$, if $D(s)=1$. For a domain $R$ with $K =\Qt{R}$, $D_K$ will denote the extension $S^{-1}D: S^{-1}A \longrightarrow S^{-1}A$ where $S = R\backslash \{ 0 \}$. $D$ is called \textit{locally nilpotent $R$-derivation ($R$-LND)} if for each $x \in A$ there exists $n \in \mathbb{N}$ such that $D^n(x) = 0$. It is well known (\textit{slice theorem}) that if $R$ is a ring containing $\mathbb{Q}$ and $D: A \longrightarrow A$ an $R$-LND having a slice $s \in A$, then $A = \KerD[s] = \KerD^{[1]}$ (see (\cite{Wright_On-The-Jacobian-Conjecture}, Proposition 2.1)); and conversely, if $D$ is irreducible and $A = \KerD^{[1]}$, then $D$ has a slice (follows from (\cite{Freudenburg-BookNew}, Principle 8(c)).	
	\medskip
	
	Let $A = R^{[n]}$ and $D: A \longrightarrow A$ an $R$-LND. The rank of $D$, denoted by $\Rk{D}$, is defined to be the least non-negative integer $r$ such that there exists a coordinate system $(X_1, X_2, \cdots, X_n)$ of $A$ satisfying $X_1, X_2, \cdots, X_{n-r} \in \KerD$. A rank-$r$ $R$-LND $D$ of $A$ is called rigid if, for any two coordinate systems $(X_1, X_2, \cdots, X_n)$ and $(Y_1, Y_2, \cdots, Y_n)$ of $A$ satisfying $X_1, X_2, \cdots, X_{n-r}, Y_1, Y_2, \cdots, Y_{n-r} \in \KerD$, we have $R[X_1, X_2, \cdots, X_{n-r}] = R[Y_1, Y_2, \cdots, Y_{n-r}]$.
	
	\medskip
	
	Affine fibrations are important objects of study in the field of affine algebraic geometry (\cite{Miyanishi_recent_dev}). A central problem in this area, due to Dolga\v{c}ev-Ve\u{\i}sfe\u{\i}ler (\cite{Dolgachev_Unipotent}), asks whether an affine fibration over a reasonably nice ring (say regular local ring) containing $\mathbb{Q}$ is a polynomial algebra. The best results in the area which directly address the above problem are the ground breaking works of Sathaye (\cite{Sat_Pol-two-var-DVR}) and Asanuma (\cite{Asanuma_fibre_ring}), which help us understand the structure of affine fibrations to a reasonable extent. Affine fibrations have close relationships with other exciting problems in affine algebraic geometry like the general epimorphism problem (\cite{Dutta-Gupta_Epimorphism-Survey}, Section 4). Indeed, investigations on affine fibrations have played a crucial role in recent breakthroughs in the Zariski cancellation problem (\cite{Neena_Threefold}; \cite{Dutta-Gupta_Epimorphism-Survey}, Section 3).
	
	\medskip
	
	As in the case of polynomial algebras over a domain $R$ containing $\mathbb{Q}$, the study of an affine fibration $A$ over $R$ may be undertaken through the studies of the LNDs of it. The questions arise: whether the kernel $B$ of an LND on $A$ is necessarily an affine fibration over the base ring $R$ and whether the given affine fibration $A$ over $R$ is also an affine fibration over the kernel $B$. When the affine fibration $A$ is a polynomial algebra, the concept of the rank of an LND is involved in some of the major results, e.g., \cite{Freudenburg_Triang_Add-Gp-Action} \cite{Daigle_NSC-Triangulability}, \cite{Daigle-Freudenburg_UFD-LND-Rank-2} (also see \cite{Freudenburg-BookNew}). However, the concept of rank has been defined only for polynomial algebras, and therefore, it is natural to ask whether a suitable notion of rank of LNDs of affine fibrations can be defined, which is consistent with the existing concept of rank of LNDs of polynomial algebras.
	
	\medskip
	
	The following two results show that when the rank of an LND of a polynomial algebra is at most two, the LND satisfies some nice properties.
	
	\begin{thm} \label{Rk2-LND_Properties_Pol-Alg}
		Let $R$ be a domain containing $\mbbQ$, $A = R^{[n]}$ and $D: A \longrightarrow A$ an $R$-LND. Then, the following hold.
		
		\begin{enumerate}
			\item [\rm (I)] Suppose that rank of $D$ is one. Then, $\KerD = R^{[n-1]}$ and $A = \KerD^{[1]}$.
			\item [\rm (II)] Suppose that rank of $D$ is two.
			
			\begin{enumerate}
				\item If $R$ is a HCF domain or a UFD, then $\KerD = R^{[n-1]}$.
				\item If $D$ is fixed point free, then $\KerD = R^{[n-1]}$ and $A = \KerD^{[1]}$.
			\end{enumerate}
		\end{enumerate}
	\end{thm}
	
	Theorem \ref{Rk2-LND_Properties_Pol-Alg}(I) follows from the property that kernel of an LND of a domain $B$ is an inert subring of $B$ and the trancendence degree of $B$ over the kernel is one. Theorem \ref{Rk2-LND_Properties_Pol-Alg}(II)(a) holds due to (\cite{AEH_Coff}, Proposition 4.1 \& Proposition 4.8); and Theorem \ref{Rk2-LND_Properties_Pol-Alg}(II)(b) appears as a corollary of (\cite{Essen_Around-Cancellation},  Remark 3.2).	
	
\medskip
	
	As a consequence of Theorem \ref{Rk2-LND_Properties_Pol-Alg} we have the following characterization of the $R$-LNDs of $R^{[3]}$ having a slice when $R$ is a PID.
	
	\begin{cor} \label{Cor_Rk2-FPF-LND_Pol-Alg}
		Let $R$ be a PID containing $\mathbb{Q}$, $A = R^{[3]}$ and $D$ a fixed point free $R$-LND of $A$. Then, the following are equivalent.
		\begin{enumerate}
			\item [\rm (I)] The rank of $D$ is at most two.
			\item [\rm (II)] $\KerD = R^{[2]}$ and $A = \KerD^{[1]}$.
			\item [\rm (III)] $D$ has a slice.
		\end{enumerate}
	\end{cor}
	
	The equivalence of (I) and (II) of Corollary \ref{Cor_Rk2-FPF-LND_Pol-Alg} follows from Theorem \ref{Rk2-LND_Properties_Pol-Alg}. (II) implies (III)  follows from the converse of \textit{slice theorem}; and (III) implies (II) follows from (\cite{Freudenburg_Rxyz-slice}, Theorem 1.1), (\cite{Sat_Pol-two-var-DVR}, Theorem 1) and (\cite{BCR_Local-Poly}, Theorem 4.4).

\begin{rem}
	\normalfont
	\begin{enumerate}
		\item [\rm (I)] It is to be noted that in Theorem \ref{Rk2-LND_Properties_Pol-Alg}, if rank of $D$ is three, then $\KerD$ need not be a polynomial ring even when $R$ is a PID and $D$ is fixed point free (see \cite{Winkelmann_free-holomorphic} or \cite{Freudenburg-BookNew}, pp.104 -- 105).
		
		\item [\rm (II)] In Corollary \ref{Cor_Rk2-FPF-LND_Pol-Alg}, if the rank of $D$ is three, then $D$ can not have slice. See (\cite{Winkelmann_free-holomorphic} or \cite{Freudenburg-BookNew}, pp.104 -- 105) for example of such LNDs. Note that Corollary \ref{Cor_Rk2-FPF-LND_Pol-Alg} holds even over one dimensional Noetherian domains containing $\mathbb{Q}$; see Lemma \ref{Lem_One-dim-Noe_Slice-implies-RankOne} for details.
	\end{enumerate}
\end{rem}
	
In section \ref{Sec_Rank}, we define residual rank and residual-variable rank of LNDs of affine fibrations. We observe that if an affine fibration is a polynomial algebra, then the rank of an LND equals to its residual rank and residual-variable rank under certain conditions (see Remark \ref{Rem_defn-Rk}(\ref{Rem_Rel_Rk-ResRk-ResVarRk})); otherwise, in general, residual rank is dominated by residual-variable rank which is dominated by rank. Further, we get results analogous to the existing results on the rank of LNDs of polynomial rings, specifically (see Corollary \ref{Cor_FibRk1-2LND} and Corollary \ref{Cor_FibRk_FPF-LND}).
	
	\medskip
	
	\noindent
	{\bf Theorem A:} Let $A$ be an $\A{n}$-fibration over a Noetherian domain $R$ containing $\mathbb{Q}$ and $D: A \longrightarrow A$ an $R$-LND. Then, the following hold.
	\begin{enumerate}
		\item [\rm (I)] 
		\begin{enumerate}
			\item If the residual rank of $D$ is one, then $\KerD$ is an $\A{n-1}$-fibration over $R$ and $A$ is an $\A{1}$-fibration over $\KerD$. Further, if $R$ is a UFD, then $A = \KerD^{[1]}$.
			
			\item If the residual-variable rank of $D$ is one, then $\KerD = R^{[n-1]}$ and $A$ is an $\A{1}$-fibration over $\KerD$. Further, if either $R$ is a UFD or $A$ is stably polynomial over $R$, then $A = \KerD^{[1]}$.
		\end{enumerate}
		\item [\rm (II)] 
		\begin{enumerate}
			\item If the residual rank of $D$ is two and $R$ is a UFD, then $\KerD = B^{[1]}$ for some $\A{n-2}$-fibration $B$ over $R$.
			\item If the residual-variable rank of $D$ is two and $R$ is a UFD, then $\KerD = R^{[n-1]}$.
		\end{enumerate}
		
		\item [\rm (III)] Suppose that $D$ is fixed point free and the residual rank of $D$ is at most two, then $D$ has a slice.
	\end{enumerate}
	
	From the above result which is parallel to Theorem \ref{Rk2-LND_Properties_Pol-Alg}, we see that triviality of affine fibration is same as having LNDs with certain residual-variable rank. As an immediate application of Theorem A, we get a characterization of the LNDs of $\A{3}$-fibrations with slice as follows (see Corollary \ref{Cor_FPF_LND_on_A^3}). The result is analogous to Corollary \ref{Cor_Rk2-FPF-LND_Pol-Alg}.	
	\medskip
	
	\noindent
	\textbf{Corollary B:} Let $R$ be a Noetherian domain containing $\mathbb{Q}$, $A$ an $\A{3}$-fibration over $R$ and $D: A \longrightarrow A$ a fixed point free $R$-LND. Then, the following are equivalent.
	\begin{enumerate}
		\item [\rm (I)] The residual rank of $D$ is at most two.
		\item [\rm (II)] $\KerD$ is an $\A{2}$-fibration over $R$ and $A$ is an $\A{1}$-fibration over $\KerD$.
		\item [\rm (III)] $D$ has a slice.
	\end{enumerate}
	Further, if the residual rank of $D$ is three, then $\KerD$ need not be an $\A{2}$-fibration over $R$ (see Example \ref{Ex_Winkelmann-1}).
		
	\medskip	

	Rigidity of LNDs of polynomial algebras plays an important role in the study of triangulability of derivations using its rank. Keshari and Lokhande proved the following result on rigidity of LNDs (\cite{Swapnil_Rigidity}, Theorem 3.1 \& Corollary 3.2) as an extension of a result by Daigle (\cite{Daigle_NSC-Triangulability}, Corollary 3.4 \& Theorem 2.5).
	\begin{thm} \label{Swapnil_Rigidity}
		Let $R$ be a domain containing $\mathbb{Q}$ with $\Qt{R} = K$ and $D$ an $R$-LND of $R^{[n]}$ such that the rank of $D$ equals to the rank of $D_K$. If $D_K$ is rigid, then $D$ is also rigid. Consequently, an $R$-LND $D$ of $R^{[3]}$ is rigid if the rank of $D$ equals to the rank of $D_K$.
	\end{thm}
	
	In section \ref{Sec_Rigidity}, we define residual-rigidity using residual rank and show that our notion of rigidity of LNDs of affine fibrations also enjoys similar property as in Theorem \ref{Swapnil_Rigidity}, specifically (see Theorem \ref{Thm_Rigidity-D-and-Dk-on-AnFib} and Corollary \ref{Cor_All-Rigid-in-A3fib}).
	
	\medskip
	
	\noindent
	{\bf Theorem C:} Let $A$ be an $\A{n}$-fibration over a Noetherian domain $R$ containing $\mathbb{Q}$ with $\Qt{R} =K$ and $D: A \longrightarrow A$ an $R$-LND. If the residual rank of $D$ equals to the rank of $D_K$ and $D_K$ is rigid, then $D$ is residually rigid. Consequently, if $n=3$ and the residual rank of $D$ equals to the rank of $D_K$, then $D$ is residually rigid.
	
	\medskip
	
	In section \ref{Sec_Examples}, we discuss a few examples of LNDs of affine fibrations and calculate their residual ranks and residual-variable ranks.
	
	\section{Preliminaries} \label{Sec_Prelim}
	In this section, we fix notation, define terminologies and state some preliminary results.
	
	\subsection{\textbf{Notation:}}
	
	For a ring $R$, an $R$-module $M$ and an $R$-algebra $A$,\\
	
	$
	\begin{array}{lll}
	\Sym{R}{M} & : & \text{Symmetric algebra of $M$ over $R$}.\\
	\Omega_R(A) & : & \text{Universal module of $R$-differentials of $A$}.\\
	\trdeg{A}{R} & : & \text{Transcendence degree of $A$ over $R$, where $R \subseteq A$ are domains}.\\
	\end{array}
	$
	
	\subsection{\textbf{Definitions:}}
	\begin{enumerate}
		\item [\rm (1)] A domain $R$ is called a \textit{HCF domain} if, for any two elements $a,b$ in $R$, the ideal $(a) \cap (b)$ is principal. \textit{HCF domains} are often called \textit{GCD domains}.
		
		\item [\rm (2)] A subring $R$ of a ring $A$ is called a \textit{retract} of $A$, if there exists a ring homomorphism $\phi: A \longrightarrow R$ such that $\phi(r) = r$ for all $r \in R$. 
		
		\item [\rm (3)] A subring $B$ of a domain $A$  is called \textit{inert} in $A$, if $fg \in B$ implies $f, g \in B$ for all $f,g \in A \backslash \{ 0 \}$. 
		
		\item [\rm (4)] Let $A$ be an $\A{n}$-fibration over a ring $R$. An $m$-tuple of elements $\ul{W} := (W_1, W_2, \cdots, W_m)$ in $A$ is called an $m$-tuple \textit{residual variable} of $A$ if they are algebraically independent over $R$, and $A \otimes_R k(P) = (R[\ul{W}] \otimes_R k(P))^{[n-m]}$ for all $P \in \Spec{R}$.
		
	\end{enumerate}
	\subsection{\textbf{Preliminary results:}}
	First, we observe the following result.
	
	\begin{lem} \label{Lem_One-dim-Noe_Slice-implies-RankOne}
		
		Let $R$ be a one dimensional Noetherian domain containing $\mathbb{Q}$, $A = R^{[3]}$ and $D$ a fixed point free $R$-LND of $A$. Then, the following are equivalent.
		\begin{enumerate}
			\item [\rm (I)] $\KerD = R^{[2]}$ and $A = \KerD^{[1]}$.
			\item [\rm (II)] $D$ has a slice.
		\end{enumerate}
	\end{lem}
	\begin{proof}
		\ul{(I) $\implies$ (II):} Follows from the converse of \textit{slice theorem}.
		
		\medskip
		\ul{(II) $\implies$ (I):} Since $D$ has a slice, by the \textit{slice theorem} we have $A = \KerD^{[1]}$, and therefore, by (\cite{Freudenburg_Rxyz-slice}, Theorem 1.1) it follows that $\KerD$ is an $\A{2}$-fibration over $R$. Since $R$ is a one dimensional Noetherian domain containing $\mathbb{Q}$, by (\cite{Asan-Bhatw_Struct-A2-fib}, Theorem 3.8) there exists $W \in \KerD$ such that $\KerD$ is an $
		\A{1}$-fibration over $R[W]$, and hence, by (\cite{DD_residual}, Corollary 3.18) $W$ is a residual variable of $\KerD$. Since $\KerD^{[1]} = R^{[3]}$, by (\cite{DD_residual}, Lemma 2.1 \& Corollary 3.19) we get $\KerD = R[W]^{[1]} = R^{[2]}$. This completes the proof.
	\end{proof}	

	We now list down some properties of inert subrings and retracts.
	\begin{lem} \label{Lem_inert-rings-and-retracts} Let $B \subseteq A$ be domains.
		\begin{enumerate}
			\item [\rm (i)] If $B$ is inert in $A$ and $C$ is such that $B \subseteq C \subseteq A$, then $B$ is inert in $C$.
			\item [\rm (ii)] If $B$ is inert in $A$, then $B$ is algebraically closed in $A$.
			\item [\rm (iii)] Let $B_1 \subseteq B_2 \subseteq A$ be domains such that $B_1$ is inert in $A$. Then, $B_1$ is inert in $B_2$. If $\trdeg{A}{B_1} = \trdeg{A}{B_2} < \infty$, then $B_1 = B_2$.
			\item [\rm (iv)] An inert subring of a HCF domain (UFD) is a HCF domain (UFD); and a polynomial algebra over a HCF domain (UFD) is a HCF domain (UFD).
			\item [\rm (v)] Retract of a UFD is a UFD.
		\end{enumerate}
	\end{lem}
	\begin{proof}
		Proofs of (i), (ii), (iii) and (iv) are easy. For a proof of (v), one may refer to the arguments by Ed Enochs mentioned in (\cite{Eakin-Heinzer_A-Cncl-prob}, p.69).
	\end{proof}

Next, we state a few properties of affine fibrations.
	
	\begin{lem} \label{Lem_AffFib-FF}
		Let $R$ be a ring and $A$ an $\A{n}$-fibration over $R$, then $A$ is faithfully flat over $R$.
	\end{lem}
	\begin{proof}
		Since $A$ is flat over $R$ and $A \otimes_R R/\mathfrak{m} = A \otimes_R k(\mathfrak{m}) = k(\mathfrak{m})^{[n]} \ne (0)$ for each maximal ideal $\mathfrak{m}$ of $R$, it follows that $A$ is faithfully flat over $R$.
	\end{proof}
	
	\begin{lem} \label{Lem_AffFib-domain-Inert}
		Let $R$ be a domain and $A$ an $\A{n}$-fibration over $R$. Then, $A$ is a domain and $R$ is inert in $A$.
	\end{lem}
	
	\begin{proof}
		Since $R$ is a domain, we have $A \hookrightarrow A \otimes_R \Qt{R} = A \otimes_R k(0) = k(0)^{[n]} = \Qt{R}^{[n]}$. This shows that $A$ is a domain and $\Qt{R}$ is inert in $A\otimes_R \Qt{R}$. Let $f, g \in A \backslash \{ 0 \}$ such that $fg \in R$. We shall show that $f,g \in R$. Due to inertness of $\Qt{R}$ in $A\otimes_R \Qt{R}$, we see that $f, g \in \Qt{R}$, and therefore, $f,g \in A \cap \Qt{R}$. Suppose, $f = r / s$ for some $r \in R$ and $s \in R \backslash \{ 0 \}$, and hence, $sf = r \in R$. Since $A$ is an $\A{n}$-fibration over $R$, by Lemma \ref{Lem_AffFib-FF} it follows that $A$ is faithfully flat over $R$, and therefore, $s A \cap R = s R$. This shows that $r = sf \in sR$, and hence, $f \in R$. Similarly, we have $g \in R$. 
	\end{proof}
	
	Asanuma, in his structure theorem of affine fibrations (\cite{Asanuma_fibre_ring}, Theorem 3.4), established that for an affine fibration $A$ over a Noetherian ring $R$, the module of differentials $\Omega_R (A)$ is a projective $A$-module and $A$ can be viewed as an $R$-subalgebra of a polynomial algebra $B$ over $R$ in such a way that $A \otimes_R B$ is a symmetric $B$-algebra of the extended projective $B$-module $\Omega_R(A) \otimes_A B$. As a consequence of Asanuma's result the following can be observed.
	
	\begin{lem} \label{Lem_AffFib-Retract}
		Let $R$ be a Noetherian ring and $A$ an $\A{n}$-fibration over $R$. Then, $R$ is a retract of $A$ and $A$ is a retract of $R^{[t]}$ for some $t \in \mathbb{N}$.
	\end{lem}
	\begin{proof}
		Since $A$ is an $\A{n}$-fibration over a Noetherian ring $R$, by (\cite{Asanuma_fibre_ring}, Theorem 3.4), $\Omega_R(A)$ is a projective $A$-module and there exists $m \in \mathbb{N}$ such that $A$ is a $R$-subalgebra of $R^{[m]}$ with the property $A^{[m]} = \Sym{R^{[m]}}{\Omega_R(A) \otimes_A R^{[m]}}$. Using Lemma \ref{Lem_AffFib-FF} we see that $A$ is a faithfully flat $R$-algebra, and therefore, $R$ can be seen as a subring of $A$. Clearly, $R$ is a retract of $R^{[m]}$, and hence, $R$ is a retract of $A$. Now, we shall show that $A$ is a retract of $R^{[t]}$ for some $t \in \mathbb{N}$.
		
		\smallskip
		
		Since $\Omega_R(A)$ is a projective $A$-module, $\Omega_R(A) \otimes_A R^{[m]}$ is a projective $R^{[m]}$-module, and therefore, we have $N \oplus (\Omega_R(A) \otimes_A R^{[m]}) = (R^{[m]})^{\ell}$ for some projective $R^{[m]}$-module $N$ and $\ell \in \mathbb{N}$. From this we get
		
		$$
		\begin{array}{ll}
		R^{[m + \ell]} & = \Sym{R^{[m]}}{(R^{[m]})^\ell}\\
		~ & = \Sym{R^{[m]}}{N \oplus (\Omega_R(A) \otimes_A R^{[m]})}\\
		~ & = \Sym{R^{[m]}}{N} \otimes_{R^{[m]}} \Sym{R^{[m]}}{\Omega_R(A) \otimes_A R^{[m]}}\\
		~ & = \Sym{R^{[m]}}{N} \otimes_{R^{[m]}} A^{[m]}.
		\end{array}
		$$
		
		Since any symmetric algebra has a natural retraction to its base ring, we see that $R^{[m]}$ is a retract of $\Sym{R^{[m]}}{N}$, and therefore, $R^{[m]} \otimes_{R^{[m]}} A^{[m]} = A^{[m]}$ is a retract of $\Sym{R^{[m]}}{N} \otimes_{R^{[m]}} A^{[m]} = R^{[m + \ell]}$. Again, since $A$ is a retract of $A^{[m]}$, we see that $A$ is a retract of $R^{[m + \ell]}$.
	\end{proof}

	\section{Rank of LNDs of affine fibrations} \label{Sec_Rank}
	
	We first define rank of an LND of an affine fibration.
	
	\begin{defn} \label{Def_ResSys-ResRk}
		Let $R$ be a ring and $A$ an $\A{n}$-fibration over $R$. 
		\begin{enumerate}
			\item For an $R$-subalgebra $B$ of $A$, the sequence $(R, B, A)$ is called an $(n, r)$-residual system if $B$ is an $\A{n-r}$-fibration over $R$ and $A \otimes_R k(P) = (B \otimes_R k(P))^{[r]}$ for all prime ideals $P$ of $R$.
			
			\item Let $D : A \longrightarrow A$ be an $R$-LND. 
			
			\begin{enumerate}
				\item [\rm (i)] $D$ is said to have \textit{residual rank} $r$ if $r$ is the least non-negative integer for which there exists an $(n, r)$-residual system $(R, B, A)$ such that $B \subseteq \Ker{D}$.
				
				\item [\rm (ii)] $D$ is said to have \textit{residual-variable rank} $r$ if $r$ is the least non-negative integer for which there exists an $(n, r)$-residual system $(R, B, A)$ such that $B = R^{[n-r]}$ and $B \subseteq \Ker{D}$.
				
			\end{enumerate} 
		\end{enumerate}
		
		The residual rank and residual-variable rank of $D$ shall be denoted by $\ResRk{D}$ and $\ResVarRk{D}$, respectively. Note that $\ResRk{D}$ and $\ResVarRk{D}$ belong to $\{0, 1, \cdots, n\}$.
	\end{defn}

\begin{rem} \label{Rem_Asanuma-Bhatwadekar_ResSys-A2-fib}
	\normalfont
	Given a non-trivial $\A{n}$-fibration $A$ over a ring $R$, there may not exists an $(n, r)$-residual system $(R,B,A)$ where $1 \le r <n$ even for the case $n=2$ (see Example \ref{Ex_Hochster}). However, Asanuma and Bhatwadekar proved that (\cite{Asan-Bhatw_Struct-A2-fib}, Theorem 3.8) when $R$ is a one-dimensional Noetherian domain containing $\mathbb{Q}$ and $A$ is an $\A{2}$-fibration over $R$, then there exists $W \in A$ such that $A$ is an $\A{1}$-fibration over $R[W]$, and therefore, by Lemma \ref{Lem_Fib-to-ResSys} it follows that $(R,R[W],A)$ is a $(2, 1)$-residual system. From their result it also follows that $A$ has an $R$-LND $D$ such that $\KerD = R[W]$, and therefore, the residual rank as well as the residual variable rank of $D$ is one.
\end{rem}
	
	In view of Definition \ref{Def_ResSys-ResRk} a result on residual variables by Das and Dutta (\cite{DD_residual}, Corollary 3.6, Theorem 3.13, Corollary 3.19 \& Appendix A) can be stated as follows.
	\begin{rem} \label{DD_residual}
		Let $R$ be a Noetherian ring, $A$ an $\A{n}$-fibration over $R$ and $(R, B , A)$ an $(n, n-r)$-residual system. Then, $A$ is an $\A{n-r}$-fibration over $B$ and $\Omega_R(A) = \Omega_B(A) \oplus (\Omega_R(B) \otimes_B A)$. Further, suppose $B = R[\ul{W}] = R^{[r]}$, i.e., $\ul{W}$ is an $r$-tuple residual variable of $A$, and $\Omega_R(A)$ is a stably free $A$-module. Then, 
		\begin{enumerate}
			\item [\rm (I)] $A^{[\ell]} = B^{[n-r+\ell]}$ for some $\ell \in \mathbb{N}$.
			\item [\rm (II)] $A = B^{[1]}$, provided $n-r =1$ and $\mbbQ \hookrightarrow R$.
		\end{enumerate}
	\end{rem}
	
	It is to be noted that though Das and Dutta, in \cite{DD_residual}, proved Remark \ref{DD_residual}(II) (see \cite{DD_residual}, Corollary 3.19) over Noetherian domains containing $\mathbb{Q}$, from their proof it follows that Remark \ref{DD_residual}(II) holds over Noetherian rings (not necessarily domains) containing $\mathbb{Q}$.

	\begin{rem} \label{Rem_defn-Rk} \normalfont 
		Let $R$ be a ring, $A$ an $\A{n}$-fibration over $R$ and $D: A \longrightarrow A$ an $R$-LND. Then, the following can be observed from Definition \ref{Def_ResSys-ResRk}.
		\begin{enumerate}
			\item \label{Rem-ResVarRk-by-ResVar}\textbf{Residual-variable rank is determined by residual variables:} From the definition of residual-variable rank and residual variable, it directly follows that $D$ has residual-variable rank $r$ if and only if $A$ has an $(n-r)$-tuple residual variable $(W_1, W_2, \cdots, W_{n-r})$ over $R$ satisfying $W_1, W_2, \cdots, W_{n-r} \in \KerD$.
			
			\item \label{Rem_ResSys-tower-AffFib_Inert-Chain} \textbf{Residual system implies tower of affine fibrations:} If $R$ is Noetherian and $(R, B, A)$ an $(n, r)$-residual system, then by Remark \ref{DD_residual} we see that $A$ is an $\A{r}$-fibration over $B$. If we further assume that $R$ is a domain, then by Lemma \ref{Lem_AffFib-domain-Inert} it follows that $R$ is inert in both $B$ and $A$, and $B$ is inert in $A$. 
			
			\item \label{Rem_ResRk-ResVarRk-zero}
			\textbf{Condition for residual rank and residual-variable rank of an LND to be zero:} If $R$ is a domain, then it is easy to see from Lemma \ref{Lem_AffFib-domain-Inert} and Lemma \ref{Lem_inert-rings-and-retracts} that $\ResRk{D} =0$ if and only if $D =0$; and $\ResVarRk{D} =0$ if and only if $D =0$ and $A$ has an $n$-tuple residual variable.
			
			\item \label{Rem_Rel_Rk-ResRk-ResVarRk}\textbf{Relation between residual rank, residual-variable rank and rank:} If $R$ is a domain with $K = \Qt{R}$ and $D$ is non-trivial, then the following hold.
			
			\begin{enumerate}
				\item Clearly, $\ResRk{D} \le \ResVarRk{D}$. Now, suppose that $\ResVarRk{D} = 1$. Since $D$ is non-trivial, $\ResRk{D} \ne 0$, and therefore, $0 < \ResRk{D} \le \ResVarRk{D} =1$, i.e., $\ResRk{D} = \ResVarRk{D} =1$. 
				
				\item Since $A \otimes_R K = K^{[n]}$, it directly follows from the definition that $\Rk{D_K} = \ResRk{D_K} = \ResVarRk{D_K} \le \ResRk{D}$. 
				
				\item Suppose $A = R^{[n]}$. Clearly, $\ResVarRk{D} \le \Rk{D}$, and therefore, $\Rk{D_K} \le \ResRk{D} \le \ResVarRk{D} \le \Rk{D}$. Hence, if we suppose that $\Rk{D} = \Rk{D_K}$, then we have $\ResRk{D} = \ResVarRk{D} = \Rk{D}$. Further, if $R$ is Noetherian and $\ResVarRk{D} = 1$, then by a residual variable result of Bhatwadekar and Dutta (\cite{BD_RES}, Remark 3.4) we get $\ResRk{D} = \ResVarRk{D} = \Rk{D} = 1$. 
			\end{enumerate}
		\end{enumerate}
	\end{rem}
	
	Now, we discuss the case residual rank and residual-variable rank at most two. Collectively, the discussion proves Theorem A. At first we observe a few results on residual systems.
	
	\begin{lem} \label{Lem_Fib-to-ResSys} Let $R$ be a domain, $A$ an $R$-algebra and $B$ an $R$-subalgebra of $A$ such that $B$ is an $\A{n-r}$-fibration over $R$ and $A$ is an $\A{r}$-fibration over $B$. Then, the following hold.
		
		\begin{enumerate}
			\item [\rm (I)] If $r =1$, then $(R, B, A)$ is an $(n, 1)$-residual system.
			\item [\rm(II)] If $R$ contains $\mathbb{Q}$, $n=3$ and $r =2$, then $(R, B, A)$ is a $(3, 2)$-residual system.
		\end{enumerate}
	\end{lem}
	
	\begin{proof} Note that $A$ is finitely generated and flat over $R$. Since $A$ is an $\A{r}$-fibration over $B$, we have $A \otimes_R k(P)$ is an $\A{r}$-fibration over $B \otimes_R k(P) = k(P)^{[n-r]}$ for all $P \in \Spec{R}$. This implies that $\trdeg{A \otimes_R k(P)}{B \otimes_R k(P)} =r$ for all $P \in \Spec{R}$.
		
		\smallskip
		
		(I): Let $r=1$. Fix $P \in \Spec{R}$. By Lemma \ref{Lem_AffFib-Retract} and Lemma \ref{Lem_inert-rings-and-retracts} we see that $B\otimes_R k(P) \subseteq A \otimes_R k(P) \subseteq B\otimes_R k(P)^{[t]}$ is a chain of UFDs for some $t \in \mathbb{N}$. Since $\trdeg{A \otimes_R k(P)}{B \otimes_R k(P)} =1$, a result of Abhyankar et al. (\cite{AEH_Coff}, Proposition 4.1) implies $A \otimes_R k(P) = (B \otimes_R k(P))^{[1]} = k(P)^{[n]}$. This shows that $A$ is an $\A{n}$-fibration over $R$, and hence, $(R, B, A)$ is an $(n, 1)$-residual system.
		
		\medskip
		
		(II): Assume that $R \hookleftarrow \mathbb{Q}$, $n=3$ and $r=2$. Fix $P \in \Spec{R}$. Since $B \otimes_R k(P) = k(P)^{[1]}$ is a PID, by a result of Sathaye (\cite{Sat_Pol-two-var-DVR}, Theorem 1) and a result of Bass et al. (\cite{BCR_Local-Poly}, Theorem 4.4) it follows that $A \otimes_R k(P) = (B \otimes_R k(P))^{[2]} = k(P)^{[3]}$. This shows that $A$ is an $\A{3}$-fibration over $R$, and hence, $(R, B, A)$ is a $(3, 2)$-residual system.
	\end{proof}
	
	\begin{thm} \label{Thm_Residual-Sys-1-2}
		Let $A$ be an $\mathbb{A}^n$-fibration over a domain $R$ containing $\mathbb{Q}$, $D: A \longrightarrow A$ a non-trivial $R$-LND and $(R,B,A)$ be an $(n, r)$-residual system such that $B \subseteq \KerD$. Then, the following hold.
		\begin{enumerate}
			\item [\rm (I)] If $r=1$ and $R$ is Noetherian, then $A$ is an $\A{1}$-fibration over $B$ and $\KerD = B$, i.e., $A$ is an $\A{1}$-fibration over $\KerD$ and $\KerD$ is an $\A{n-1}$-fibration over $R$. If we further assume that either $A$ is stably polynomial over $R$ and $B = R^{[n-1]}$ or $D$ is fixed point free, then $A = \KerD^{[1]}$.
			
			\item [\rm (II)] If $r=2$, $R$ is Noetherian and $D$ is fixed point free, then $A = \KerD^{[1]}$ and $\KerD$ is an $\A{1}$-fibration over $B$ as well as an $\A{n-1}$-fibration over $R$. If we further assume that  $A$ is stably polynomial over $B$, then $\KerD = B^{[1]}$.     
			\item [\rm (III)] Suppose that $R$ is a Noetherian UFD.
			\begin{enumerate}
				\item If $r=1$, then $A = \KerD^{[1]}$ and $\KerD = B$. 
				\item If $r=2$, then $A$ is an $\A{2}$-fibration over $B$ and $\KerD = B^{[1]}$.
			\end{enumerate}
			\item [\rm (IV)] Suppose that $R$ is a HCF domain and $A$ is stably polynomial over $R$ as well as over $B$.
			\begin{enumerate}
				\item If $r=1$, then $A = \KerD^{[1]}$ and $\KerD = B$.
				
				\item If $r=2$, then $\KerD =B^{[1]}$.
			\end{enumerate}
		\end{enumerate}
	\end{thm}
	
	\begin{proof} Since $D \ne 0$ and $A$ is a domain, by (\cite{Freudenburg-BookNew}, Principle 1 and Principle 11), $\trdeg{A}{\KerD} = 1$ and $\KerD$ is inert in $A$.
		
		\smallskip
		
		(I): Let $R$ be Noetherian and $r=1$. By Remark \ref{Rem_defn-Rk}(\ref{Rem_ResSys-tower-AffFib_Inert-Chain}) $A$ is an $\A{1}$-fibration over $B$ and $B$ is inert in $A$. Since $B \subseteq \KerD \subseteq A$ and $\trdeg{A}{B} = \trdeg{A}{\KerD} = 1$, using Lemma \ref{Lem_inert-rings-and-retracts} we get $\KerD = B$. 
		
		\smallskip
		
		Now we further assume that $A$ is stably polynomial over $R$. By (\cite{DD_residual}, Lemma 2.1) $\Omega_R(A)$ is stably free over $A$. Since $(R, B, A)$ is a residual system and $\KerD = B = R^{[n-1]}$, applying Remark \ref{DD_residual} we get $A = \KerD^{[1]}$.
		
		\smallskip
		
		Again, along with the hypotheses $r=1$ and $R$ is Noetherian, if we further suppose that $D$ is fixed point free, then by a result of Kahoui-Ouali (\cite{Kahoui_Cancellation}, Corollary 2.5) we get that $A = \KerD^{[1]}$.
		
		\medskip
		
		(II): Let $R$ be Noetherian, $r=2$, and $D$ a fixed point free $R$-LND. Since $(R,B,A)$ is an $(n, 2)$-residual system, from Remark \ref{Rem_defn-Rk}(\ref{Rem_ResSys-tower-AffFib_Inert-Chain}) we have $A$ is an $\A{2}$-fibration over $B$. Now, since $B \subseteq \KerD$, we see that $D$ is a $B$-LND of $A$, and hence by a result of Babu and Das (\cite{Das-Raman_Struct_A2-fib_FPF-LND}, Remark 4.9) it follows that $A = \KerD^{[1]}$ and $\KerD$ is an $\A{1}$-fibration over $B$. Since $B$ is an $\A{n-2}$-fibration over $R$, by Lemma \ref{Lem_Fib-to-ResSys}(I) we get $(R,B, \KerD)$ is an $(n-1, 1)$-residual system and which implies that $\KerD$ is an $\A{n-1}$-fibration over $R$. If we further assume that $A$ is stably polynomial over $B$, then by (\cite{Das-Raman_Struct_A2-fib_FPF-LND}, Remark 4.9) $\KerD = B^{[1]}$.
		
		\medskip
		
		(III): Let us assume that the hypothesis holds. Since both $A$ and $B$ are affine fibrations over $R$, by Lemma \ref{Lem_AffFib-Retract} and Lemma \ref{Lem_inert-rings-and-retracts} we see that both $A$ and $B$ are UFDs. Again, since $\KerD$ is inert in $A$, by Lemma \ref{Lem_inert-rings-and-retracts} it follows that $\KerD$ is also a UFD.
		
		\begin{enumerate}
			\item [\rm (a)] Let $r=1$. By (I) we have $\KerD = B$. Since $A$ is an $\A{1}$-fibration over $\KerD$, by (\cite{Asanuma_fibre_ring}, Theorem 3.4) we find a $t \in \mathbb{N}$ such that $A$ is a $\KerD$-subalgebra of $\KerD^{[t]}$. This shows that $\KerD \subseteq A \subseteq \KerD^{[t]}$ is a chain of UFDs where $\trdeg{A}{\KerD} = 1$. Therefore, by (\cite{AEH_Coff}, Proposition 4.1) we conclude that $A = B^{[1]}$.
			
			\item [\rm (b)] Let $r=2$. By Remark \ref{Rem_defn-Rk}(\ref{Rem_ResSys-tower-AffFib_Inert-Chain}), $A$ is an $\A{2}$-fibration over $B$ and therefore, using (\cite{Asanuma_fibre_ring}, Theorem 3.4) we get an $\ell \in \mathbb{N}$ such that $A$ is an $B$-subalgebra of $B^{[\ell]}$. Notice that $B \subseteq \KerD \subseteq A \subseteq B^{[\ell]}$ is a chain of UFDs where $\trdeg{\KerD}{B} =1$. Therefore, by (\cite{AEH_Coff}, Proposition 4.1), we get $\KerD = B^{[1]}$.
		\end{enumerate}
		
		(IV): We assume the hypothesis. Since $A$ is stably polynomial over both $R$ and $B$, there exist $s, t \in \mathbb{N}$ such that $A^{[s]} = R^{[n+s]}$ and $A^{[t]} = B^{[r+t]}$. Since $R \subseteq A \subseteq R^{[n+s]}$ and $B \subseteq A \subseteq B^{[r+t]}$, by Lemma \ref{Lem_inert-rings-and-retracts} both $R$ and $B$ are inert in $A$, and therefore, repeated application of Lemma \ref{Lem_inert-rings-and-retracts} implies that $A$, $B$ and $B^{[m]}$ are HCF domains for all $m \in \mathbb{N}$.
		\begin{enumerate}
			\item [\rm (a)] Let $r =1$. By (I) we have $\KerD = B$, and hence $A^{[t]} = B^{[t+1]} = \KerD^{[t+1]}$. Clearly, $A$ is inert in $\KerD^{[t+1]}$ and $\KerD \subseteq A \subseteq \KerD^{[t+1]}$ is a chain of HCF domains. By (\cite{AEH_Coff}, Proposition 4.8) we get $A = \KerD^{[1]}$.
			
			\item [\rm (b)] Let $r=2$. Since $A$ is stably polynomial over $B$, we see that $B \subseteq \KerD \subseteq A^{[t]} = B^{[t+2]}$ is a chain of HCF domains for some there $t \in \mathbb{N}$. Note that since $\KerD$ is inert in $A$, it is also inert in $A^{[t]}$. Now, by (\cite{AEH_Coff}, Proposition 4.8) we conclude that $\KerD = B^{[1]}$.
		\end{enumerate}
	\end{proof}	
	
	As a consequence of Theorem \ref{Thm_Residual-Sys-1-2}, we get the following analogue of Theorem \ref{Rk2-LND_Properties_Pol-Alg}(I) and (II)(a).
	
	\begin{cor} \label{Cor_FibRk1-2LND}
		Let $A$ be an $\A{n}$-fibration over a Noetherian domain $R$ containing $\mathbb{Q}$ and $D: A \longrightarrow A$ a non-trivial $R$-LND. Then, the following hold.
		\begin{enumerate}
			\item [\rm (I)] Suppose that $\ResRk{D} = 1$. Then, $\KerD$ is an $\A{n-1}$-fibration over $R$ and $A$ is an $\A{1}$-fibration over $\KerD$. Furthermore, if either $A$ is assumed to be stably polynomial over $R$ or $R$ is assumed to be a UFD, then $A = \KerD^{[1]}$.
			
			\item [\rm (II)] Suppose that $\ResVarRk{D} = 1$. Then, $\KerD = R^{[n-1]}$ and $A$ is an $\A{1}$-fibration over $\KerD$. Furthermore, if either $R$ is assumed to be a UFD or $A$ is assumed to be stably polynomial over $R$, then $\KerD = R^{[n-1]}$ and $A = \KerD^{[1]}$, i.e., $A = R^{[n]}$.
			
			\item [\rm (III)] Suppose that $\ResRk{D} =2$ and $R$ is a UFD. Then, $\KerD = B^{[1]}$ for some $\A{n-2}$-fibration $B$ over $R$.

			\item [\rm (IV)] Suppose that $\ResVarRk{D} = 2$ and $R$ is a UFD. Then, $\KerD = R^{[n-1]}$.
			
		\end{enumerate}
	\end{cor}
	\begin{proof}
		(I): Since $\ResRk{D} =1$, there exists an $(n, 1)$-residual system $(R, B, A)$ such that $B \subseteq \KerD$, and therefore, the result follows due to Theorem \ref{Thm_Residual-Sys-1-2}[(I) \& (III)].
		
		\medskip
		
		(II): Since $\ResVarRk{D} =1$, there exists $B = R^{[n-1]} \subseteq \KerD$ such that $(R, B, A)$ is an $(n, 1)$-residual system, and therefore, the conclusion follows from Theorem \ref{Thm_Residual-Sys-1-2}[(I) \& (III)].
		
		\medskip
		
		(III): Since $\ResRk{D} = 2$, there exists an $(n, 2)$-residual system $(R, B, A)$ such that $B \subseteq \KerD$, and therefore, by Theorem \ref{Thm_Residual-Sys-1-2}(III) we get the result.
		
		\medskip
		
		(IV): Since $\ResVarRk{D} =2$, there exists $B = R^{[n-2]} \subseteq \KerD$ such that $(R, B, A)$ is an $(n, 2)$-residual system, and therefore, by Theorem \ref{Thm_Residual-Sys-1-2}(III) the desired result holds.
	\end{proof}
	
	As an immediate corollary of Theorem \ref{Thm_Residual-Sys-1-2}(III) we observe the following result.
	
\begin{cor} \label{Cor_Ker-Triv_A2-fib_over_UFD} Let $R$ be a Noetherian UFD containing $\mathbb{Q}$, $A$ an $\A{2}$-fibration over $R$ and $D: A \longrightarrow A$ a non-trivial $R$-LND. Then, $\KerD = R^{[1]}$. 
\end{cor}

\begin{proof} Clearly, $(R, R, A)$ is a $(2, 2)$-residual system such that $R \subseteq \KerD$, and therefore, by Theorem \ref{Thm_Residual-Sys-1-2}(III) we directly get $\KerD = R^{[1]}$.
\end{proof}
	
	Next, we prove an analogue of Theorem \ref{Rk2-LND_Properties_Pol-Alg}(II)(b).
		
	\begin{cor} \label{Cor_FibRk_FPF-LND}
		Let $A$ be an $\mathbb{A}^n$-fibration over a Noetherian domain $R$ containing $\mathbb{Q}$ and $D: A \longrightarrow A$ a fixed point free $R$-LND. Then, the following hold.
		\begin{enumerate}
			\item [\rm (I)] If $\ResRk{D} \le 2$, then $A = \KerD^{[1]}$ and $\KerD$ is an $\A{n-1}$-fibration over $R$, i.e., $\ResRk{D} = 1$.
			
			\item [\rm (II)] If $\ResVarRk{D} = 1$, then $A = \KerD^{[1]}$ and $\KerD = R^{[n-1]}$.
			
			\item [\rm (III)] If $\ResVarRk{D} = 2$, then $A = \KerD^{[1]}$, $\KerD$ is an $\A{1}$ fibration over $R^{[n-2]}$, and furthermore, if $A$ is stably polynomial over $R$, then $\KerD = R^{[n-1]}$.
		\end{enumerate}
	\end{cor}        
	
	\begin{proof}
		(I): Let $\ResRk{D} \le 2$. Then, there exists an $(n, r)$-residual system $(R,B,A)$ such that $B \subseteq \KerD$ where either $r=1$ or $r=2$. Now, from Theorem \ref{Thm_Residual-Sys-1-2}[(I) \& (II)] it follows that $A = \KerD^{[1]}$ and $\KerD$  is an $\A{n-1}$-fibration over $R$, which, by Lemma \ref{Lem_Fib-to-ResSys}, is equivalent to say that $\ResRk{D} =1$.
		
		\medskip
		
		(II): Let $\ResVarRk{D} =1$. Then, there exists an $(n, 1)$-residual system $(R,B,A)$ such that $B \subseteq \KerD$ and $B = R^{[n-1]}$, and therefore, by Theorem \ref{Thm_Residual-Sys-1-2}(I) we get $A = \KerD^{[1]}$ and $\KerD = B = R^{[n-1]}$.
		
		\medskip
		
		(III): Let $\ResVarRk{D} = 2$. Then, there exists an $(n, 2)$-residual system $(R,B,A)$ such that $B \subseteq \KerD$ and $B = R^{[n-2]}$, and therefore, the conclusion follows from Theorem \ref{Thm_Residual-Sys-1-2}(II).
	\end{proof}

	\begin{rem} \label{Rem_FPF-and-ResRk2-imply-ResRk1_and_Slice-ZCP}
		\normalfont
		\begin{enumerate}
			\item \label{Rem_FPF-and-ResRk2-imply-ResRk1}The phenomenon in Corollary \ref{Cor_FibRk_FPF-LND}(I) is very specific for fixed point free LNDs, i.e., if the LND $D$ is not fixed point free then the condition $\ResRk{D} \le 2$ need not imply $\ResRk{D} =1$. One may look at Example \ref{Ex_BD_A2} for details.
			
			\item \label{Rem_Slice-ZCP} It is to be noted that in the hypothesis of Corollary \ref{Cor_FibRk_FPF-LND} if we replace the condition ``$D$ is fixed point free'' by the stronger condition ``$D$ has a slice, i.e., $A = \KerD^{[1]}$'', then the finite generation and flatness of $A$ over $R$ will imply the finite generation and flatness of $\KerD$ over $R$, and further, we shall have $k(P)^{[n]} = A \otimes_R k(P) = (\KerD \otimes_R k(P))^{[1]}$ for all $P \in \Spec{R}$.  But, since the Zariski cancellation problem\footnote{\textbf{Zariski  cancellation problem:} Let $k$ be a field and $A$ an $n$-dimensional affine $k$-algebra such that $A^{[m]} = k^{[m+n]}$. Is then $A = k^{[n]}$?} is open in dimension $n \ge 3$ over fields containing $\mathbb{Q}$, we can not conclude that $\KerD \otimes_R k(P) = k(P)^{[n-1]}$ for all $P \in \Spec{R}$, i.e., we can not conclude that $\KerD$ is an $\A{n-1}$-fibration over $R$.  However, from the above discussion it is easy to see the following.
			
			\smallskip
			
			\textit{Let $R$ be a domain containing $\mathbb{Q}$, $A$ an $\A{n}$-fibration over $R$ and $D: A \longrightarrow A$ an $R$-LND having a slice. If the Zariski cancellation problem has affirmative solution in dimension $n$ over fields containing $\mathbb{Q}$, then, $\KerD$ is an $\A{n-1}$-fibration over $R$.}
			
		\end{enumerate}
	\end{rem}

	We now prove Corollary B.
	
	\begin{cor} \label{Cor_FPF_LND_on_A^3}
		Let $R$ be a Noetherian domain containing $\mathbb{Q}$, $A$ an $\A{3}$-fibration over $R$ and $D: A \longrightarrow A$ a fixed point free $R$-LND. Then, the following are equivalent.
		
		\begin{enumerate}
			\item [\rm (I)] $D$ has a slice.
			\item [\rm (II)] $\ResRk{D} = 1$.
			\item [\rm (III)] $\ResRk{D} \le 2$.
			\item [\rm (IV)] $\KerD$ is an $\A{2}$-fibration over $R$ and $A = \KerD^{[1]}$.
			\item [\rm (V)] $\KerD$ is an $\A{2}$-fibration over $R$ and $A$ an $\A{1}$-fibration over $\KerD$.
			\item [\rm (VI)] $\KerD$ is Noetherian and $A$ is an $\A{1}$-fibration over $\KerD$.
		\end{enumerate}
		
		Further, if $\ResRk{D} =3$, then $\KerD$ need not be an $\A{2}$-fibration over $R$.
	\end{cor}
	
	\begin{proof}
		
		\underline{(I) $\implies$ (II):} Since $D$ has a slice, and since the Zariski cancellation problem has affirmative answer in dimension two over fields containing $\mbbQ$ (follows from \cite{Miyanishi_cyllinder}, \cite{Fujita_Zariski-Problem}, and \cite{Kambayashi_K2Form}), from Remark \ref{Rem_FPF-and-ResRk2-imply-ResRk1_and_Slice-ZCP}(\ref{Rem_Slice-ZCP}) it follows that $\KerD$ is an $\A{2}$-fibration over $R$. Since $A = \KerD^{[1]}$, we see that $\ResRk{D} \le 1$. Since $D$ is non-trivial, we have $\ResRk{D} =1$.
		
		\smallskip
		
		\underline{(II) $\implies$ (III)}, \underline{(IV) $\implies$ (V)} and \underline{(V) $\implies$ (VI)}: Obvious.
		
		\smallskip
		
		\underline{(III) $\implies$ (IV)}: Directly follows from Corollary \ref{Cor_FibRk_FPF-LND}(I).
		
		\smallskip
		
		\underline{(VI) $\implies$ (I)}: Directly follows from (\cite{Kahoui_Cancellation}, Corollary 2.5) and the converse of \textit{slice theorem}.
		
		\bigskip
		
		Example \ref{Ex_Winkelmann-1} exhibits an $R$-LND $D$ such that $\ResRk{D} =3$, but $\KerD$ is not an $\A{2}$-fibration over $R$. This show that if $\ResRk{D} =3$, then $\KerD$ need not be an $\A{2}$-fibration over $R$. 
	\end{proof}
	
	\begin{rem}\textbf{(Recognizing coordinates of $R^{[3]}$)}
		Let $R$ be a Noetherian domain containing $\mathbb{Q}$ and $A = R^{[3]}$. Suppose that $F \in A$ is a residual variable of $A$. Then, the following are equivalent.
		
		\begin{enumerate}
			\item [\rm (I)] $F$ is a variable of $A$.
			\item [\rm (II)] $A$ has a fixed point free $R$-LND $D$ of residual rank at most two such that $F \in \KerD$.
			\item [\rm (III)] $D$ is an $R$-LND with a slice satisfying $F \in \KerD$.
		\end{enumerate}		
		
		\smallskip
		
		\normalfont
		We give a proof to the above statement.
		
		\medskip
		
		\underline{(I) $\implies$ (II)}: Suppose that $F$ is a variable of $A$, and therefore, there exists $G,H \in A$ such that $A = R[F, G,H]$. Consider the partial derivative $\partial_G$ on A. Then, one may see that $\partial_G$ is fixed point free and $F \in \Ker{\partial_G}$.
		
		\smallskip
		
		\underline{(II) $\implies$ (III)}: Directly follows from Corollary \ref{Cor_FPF_LND_on_A^3}.
		
		\smallskip
		
		\underline{(III) $\implies$ (I)}: Assume that (III) holds. Then, by Corollary \ref{Cor_FPF_LND_on_A^3} we have $A = \KerD^{[1]}$ and $\KerD$ is an $\A{2}$-fibration over $R$ such that $R[F] \subseteq \KerD$. Let $P \in \Spec{R}$. Since $F$ is a residual variable of $A$, we have $(R[F] \otimes_R k(P))^{[2]} = A \otimes_R k(P)$. Note that $A \otimes_R k(P) = (\KerD \otimes_R k(P))^{[1]}$, and therefore, we get $(\KerD \otimes_R k(P))^{[1]} = (R[F] \otimes_R k(P))^{[2]}$. Since $\KerD \otimes_R k(P)$ is an $R[F] \otimes_R k(P)$-algebra, by a cancellation result of Hamann (\cite{Haman_Invariance}, Theorem 2.8) it follows that $\KerD \otimes_R k(P) = (R[F] \otimes_R k(P))^{[1]}$, i.e., $F$ is a residual variable of $\KerD$. Since $\KerD^{[1]} = A = R^{[3]}$, by (\cite{DD_residual}, Lemma 2.1) we get $\Omega_R(\KerD)$ is stably free over $\KerD$, and therefore, by Remark \ref{DD_residual} we conclude that $\KerD = R[F]^{[1]}$, i.e., $A = R[F]^{[2]}$.
		
		\medskip
		
		In Section \ref{Sec_Examples}, we quote an example (see Example \ref{Ex_BD_A2}), due to Bhatwadekar and Dutta (\cite{BD_AFNFIB}), of an element $F \in A = R^{[3]}$ such that there exists a non-fixed point free $R$-LND $D: A \longrightarrow A$ having residual rank two with the property that $F \in \KerD$. But it is not known whether $A$ has a fixed point free $R$-LND of residual rank two such that $F$ belongs to its kernel.
	\end{rem}    
	
	\section{Rigidity of LNDs of affine fibrations} \label{Sec_Rigidity}
	
	First, we define rigidity of LNDs of affine fibrations.
	
	\begin{defn} \label{rigidity}
		Let $A$ be an $\A{n}$-fibration over a ring $R$ and $D: A \longrightarrow A$ an $R$-LND with residual rank $r$. We define $D$ to be residually rigid if, for any two $(n, r)$-residual systems $(R, B_1, A)$ and $(R, B_2, A)$ with  $B_1, B_2 \subseteq \text{Ker}(D)$ we have $B_1 = B_2$.
	\end{defn}
	
	\begin{rem}\normalfont Let $R$ be a ring, $A$ an $R$-algebra and $D : A \longrightarrow A$ an $R$-LND.
		\begin{enumerate}
			\item  If $A = R^{[n]}$ and $\Rk{D} = \ResRk{D}$, then one can see that the residual-rigidity of $D$ implies rigidity of $D$.
			
			\item If $R$ is a domain, $A = R^{[n]}$ and $\Rk{D} =1$, then it can be seen, due to inertness of $\KerD$, that $D$ is rigid. In the context of affine fibrations one may observe a similar phenomenon, also caused by the inertness of the kernel of the LNDs: If $R$ is a Noetherian domain, $A$ is an $\A{n}$-fibration over an $R$ and $\ResRk{D} =1$, then $D$ is residually rigid. 
		\end{enumerate}
	\end{rem}
	
	Before proving Theorem C, we note the following lemma which can be seen as an extension of an observation of Abhyankar et al. (\cite{AEH_Coff}, 1.7).
	
	\begin{lem} \label{Lem_AEH-inertness-Extn}
		Let $A$ be a domain and $B_1, B_2$ subdomains of $A$. Suppose $B_2$ is inert in $A$. If $b \in B_1$ is such that $bA \cap B_2 \ne \{ 0 \}$, then $b \in B_2$.  
	\end{lem}
	
	\begin{proof}
		Let $d \in bA \cap B_2$. Then, $d = bc \in B_2$ for some $c \in A$. Since $b,c \in A$ and $B_2$ is inert in $A$, we have $b,c \in B_2$.
	\end{proof}
	
	We now prove Theorem C.
	
	\begin{thm} \label{Thm_Rigidity-D-and-Dk-on-AnFib}
		Let $A$ be an $\A{n}$-fibration over a Noetherian domain $R$ and $D: A \longrightarrow A$ an $R$-LND. Suppose that $\ResRk{D} = \Rk{D_K}$. If $D_K$ is rigid, then $D$ is residually rigid. 
	\end{thm}
	
	\begin{proof} Let $\ResRk{D} = \Rk{D_K} =r$. Let us assume that $D_K$ is rigid. Since $K$ is a field, we have $\ResRk{D_K} = \Rk{D_K} =r$. Let $(R, B_1, A)$ and $(R, B_2, A)$ be two $(n, r)$-residual systems such that $B_1, B_2 \subseteq \KerD$. By Remark \ref{Rem_defn-Rk}(\ref{Rem_ResSys-tower-AffFib_Inert-Chain}) we get $A$ is an $\A{r}$-fibration over both $B_1$ and $B_2$, and both $B_1$ and $B_2$ are inert in $A$. We shall show that $B_1 = B_2$. 
		
		\smallskip
		
		Let $\underline{U} \in {B_1}^{n-r}$ and $\underline{V} \in {B_2}^{n-r}$ be such that $B_1 \otimes_R K =K[\ul{U}]$ and $B_2 \otimes_R K = K[\ul{V}]$. Since $A$ is an $\A{r}$-fibration over both $B_1$ and $B_2$, we have $A \otimes_R K = K[\underline{U}] ^{[r]} =  K[\underline{V}]^{[r]} = K^{[n]}$, and therefore, since $D_K$ is rigid and $\Rk{D_K} = r$, we have $K[\ul{U}] = K[\ul{V}]$, i.e., $B_1 \otimes_R K = K[\underline{V}]$. Suppose that $x \in B_1$. Since $B_1 \otimes_R K = K[\ul{U}] = K[\ul{V}] = B_2 \otimes_R K$, there exists $r \in R$ such that $rx \in B_2$. This shows that $rx \in xA \cap B_2$, and therefore, by Lemma \ref{Lem_AEH-inertness-Extn} we have $x \in B_2$. So, we get $B_1 \subseteq B_2$. Now, interchanging the roles of $B_1$ and $B_2$ we get $B_2 \subseteq B_1$. Hence, $B_1 = B_2$. This completes the proof.
	\end{proof}
	
	As a direct consequence of Theorem \ref{Thm_Rigidity-D-and-Dk-on-AnFib} and (\cite{Daigle_NSC-Triangulability}, Theorem 2.5) we get the following.
	
	\begin{cor} \label{Cor_All-Rigid-in-A3fib}
		Let $A$ be an $\A{3}$-fibration over a Noetherian domain $R$ and $D: A \longrightarrow A$ an $R$-LND such that $\ResRk{D} = \Rk{D_K}$, then $D$ is residually rigid.
	\end{cor}
	
	\section{Examples} \label{Sec_Examples}
	We now discuss a few examples. The first example involves a non-trivial $\A{2}$-fibration along with a fixed point free LND. 
	
	\begin{ex} \label{Ex_NonTriv_A2-Fib_ResRk-ResVarRk} \normalfont
		Let $R$ be a Noetherian domain containing $\mbbQ$, $A$ a non-trivial $\A{2}$-fibration over $R$ and $D: A \longrightarrow A$ a fixed point free $R$-LND. We shall show $\ResRk{D} =1$ and $\ResVarRk{D} =2$. 
		
		\smallskip
		
		By (\cite{Das-Raman_Struct_A2-fib_FPF-LND}, Remark 4.9) we have $A = \Ker{D}^{[1]}$ and $\Ker{D}$ is an $\A{1}$-fibration over $R$, and therefore, by Lemma \ref{Lem_Fib-to-ResSys} we see that $(R, \KerD, A)$ is a $(2, 1)$-residual system. This shows that $\ResRk{D} \le 1$. Since $D$ is non-trivial, we have $\ResRk{D} \ne 0$, and hence $\ResRk{D} =1$. By Remark \ref{Rem_defn-Rk}(\ref{Rem_Rel_Rk-ResRk-ResVarRk}) we get $1 = \ResRk{D} \le \ResVarRk{D}$. We claim that $\ResVarRk{D} =2$. Otherwise, by Corollary \ref{Cor_FibRk_FPF-LND} we shall get $\KerD = R^{[1]}$, i.e., $A = R^{[2]}$ which contradicts the fact that $A$ is a non-trivial $\A{2}$-fibration over $R$.
		
		\smallskip
		
		Note that by (\cite{Das-Raman_Struct_A2-fib_FPF-LND}, Remark 4.9) $A$ is has another $R$-LND $D_1$ such that $\Ker{D_1} = R[V] = R^{[1]}$ for some $V \in A$ and $A$ is an $\A{1}$-fibration over $R[V]$, and therefore, by Lemma \ref{Lem_Fib-to-ResSys} $(R, R[V], A)$ is a $(2, 1)$-residual system with $R[V] \subseteq \Ker{D_1}$. This shows that $\ResVarRk{D_1} \le 1$. Now, since $D_1 \ne 0$,  by Remark \ref{Rem_defn-Rk}(\ref{Rem_Rel_Rk-ResRk-ResVarRk}) we get $0 \ne \ResRk{D_1} \le \ResVarRk{D_1}$, and therefore, we have $\ResRk{D_1} = \ResVarRk{D_1} =1$.  
	\end{ex}

The next example is by Hochster (see \cite{Hochster_nonunique-coeff} or (\cite{Freudenburg-BookNew}, 10.1.5)).

\begin{ex} \label{Ex_Hochster}
	\normalfont
	Let $R = \mathbb{R}[X,Y,Z]/(X^2 + Y^2 + Z^2 -1) = \mathbb{R}[x,y,z]$ where $x,y,z$ denote the images of $X,Y,Z$ in $R$. Let $A = R[U,V,W]/(xU + yV + zW)$. One can see that $R$ is a Noetherian UFD, $A$ is a non-trivial $\A{2}$-fibration over $R$ and $A^{[1]} = R^{[3]}$. We claim that there does not exist $B \subseteq A$ such that $(R,B,A)$ is a $(2, 1)$-residual system. On the contrary, let $(R,B,A)$ be a $(2, 1)$-residual system. Since $R$ is a Noetherian domain and $B$ is an $\A{1}$-fibration over $R$, by Lemma \ref{Lem_AffFib-Retract} and Lemma \ref{Lem_inert-rings-and-retracts} we see that $R \subseteq B \subseteq R^{[m]}$ for some $m \in \mathbb{N}$ is a sequence of UFDs with $\trdeg{B}{R} =1$, and therefore, by (\cite{AEH_Coff}, Proposition 4.1) we have $B = R^{[1]}$. Now, since $A$ is stably polynomial over $R$, by (\cite{DD_residual}, Lemma 2.1) and Remark \ref{DD_residual} it follows that $A = R^{[2]}$ which is a contradiction to the fact that $A$ is a non-trivial $\A{2}$-fibration over $R$. This shows that for any non-trivial $R$-LND $D$ of $A$, the residual rank of $D$ is always two, and therefore, in view of Example \ref{Ex_NonTriv_A2-Fib_ResRk-ResVarRk}, $A$ does not have any fixed point free $R$-LND. In this context, one should note that by Corollary \ref{Cor_Ker-Triv_A2-fib_over_UFD} we have $\KerD = R^{[1]}$; however, $A$ can not be an $\A{1}$-fibration over $\KerD$. Further, note that the above arguments and observations hold true for any non-trivial stably polynomial $\A{2}$-fibration $A$ over a Noetherian UFD $R$ containing $\mathbb{Q}$.  
\end{ex}

	The following example is by Winkelmann (\cite{Winkelmann_free-holomorphic}, also see \cite{Freudenburg-BookNew}, pp.104 -- 105).
	
	\begin{ex} \label{Ex_Winkelmann-1} \normalfont 
		Let $R = \mbbC[X] = \mbbC^{[1]}$, $A = R[U,V,W] = R^{[3]}$ and $D: A \longrightarrow A$ be an $R$-LND defined by $D(U) = X$, $D(V) = U$ and $D(W) = U^2 - 2XV -1$. One can easily see that $D$ is fixed point free. It is known that $\KerD = R[f,g,h] \ne R^{[2]}$ where 
		$$
		\begin{array}{ll}
		f  =  U^2 - 2XV, \\
		g  =  XW + (1-f)U, \\
		Xh = g^2 -f(1-f)^2, \ i.e., \ h = XW^2 + 2(1-f) (UW +(1-f)V).\\
		\end{array}
		$$
		
		By Theorem \ref{Rk2-LND_Properties_Pol-Alg} it follows that $\Rk{D} =3$. We shall calculate $\ResRk{D}$ and $\ResVarRk{D}$. Note that $\KerD$ is not an $\A{2}$-fibration over $R$, otherwise by (\cite{Sat_Pol-two-var-DVR}, Theorem 1) and (\cite{BCR_Local-Poly}, Theorem 4.4) we get $\KerD = R^{[2]}$ which is a contradiction. Further, $D$ has no slice; otherwise by Corollary \ref{Cor_FPF_LND_on_A^3} it would follow that $\KerD$ is an $\A{2}$-fibration over $R$, which is a contradiction. Thus, $D$ is fixed point free without a slice, and therefore, by Corollary \ref{Cor_FibRk_FPF-LND}, we have $\ResRk{D} =3$. So, by Remark \ref{Rem_defn-Rk}(\ref{Rem_Rel_Rk-ResRk-ResVarRk}) we see that $3 = \ResRk{D} \le \ResVarRk{D} \le \Rk{D} =3$, i.e., $\ResRk{D} = \ResVarRk{D} = \Rk{D} =3$.
		
		\smallskip
		
		Please note that, in this case, we can also use Corollary \ref{Cor_FibRk_FPF-LND} to compute $\ResVarRk{D}$ directly.
	\end{ex} 
	
	We now consider an example by Bhatwadekar and Dutta (\cite{BD_AFNFIB}, Example 4.13) (also see \cite{Venereau_Thesis}).
	
	\begin{ex}\label{Ex_BD_A2} \normalfont 
		Let $\mathbb{F}$ be a field containing $\mbbQ$, $R = \mathbb{F}[\pi]_{(\pi)}$ and $A = R[X,Y,Z]$. Set $F := \pi^2 X + \pi Y (YZ + X + X^2) + Y$. One can check that $A \otimes_R k (P)= (R[F] \otimes_R k(P))^{[2]}$ for all $P \in \Spec{R}$, i.e., $F$ is a residual variable of $A$, and therefore, by Remark \ref{DD_residual}, $A$ is an $\A{2}$-fibration over $R[F]$, and hence, by (\cite{Asanuma_fibre_ring}, Theorem 3.4), $A$ is stably polynomial over $R[F]$. It is not known whether $A = R[F]^{[2]}$. 
		
		\smallskip
		
		Define an $R$-LND $D$ of $A$ by $D(X) = Y^2$, $D(Y) = 0$ and $D(Z) = - (\pi + Y + 2XY )$. Then, $R[F] \subseteq \Ker{D}$. We shall find $\Rk{D}$, $\ResRk{D}$ and $\ResVarRk{D}$. Clearly, $D$ is irreducible and triangular. By Corollary \ref{Cor_Ker-Triv_A2-fib_over_UFD} we get $\KerD = R[F]^{[1]} = R^{[2]}$. We now show that $D$ is not fixed point free. On the contrary, assume that $D$ is fixed point free, and therefore, there exists $f_1, f_2, f_3 \in R[X,Y,Z]$ such that $D(X)f_1 + D(Y)f_2 + D(Z) f_3 = 1$. Since $D(Y) = 0$, we have $D(X)f_1 + D(Z) f_3 = 1$, i.e., $Y^2 f_1 - (\pi + Y + 2XY)f_3 = 1$. Hence, in $A/YA = R[X,Z]$ we get $-\pi f_3 =1$, i.e., $\pi$ is a unit in $R[X,Z]$, giving a contradiction to the fact that $\pi$ is a prime in $R$.
		
		\smallskip
		
		As $A = R[X,Y,Z]$ and $Y \in \KerD$, we see that $\Rk{D} \le 2 $. Since $R$ is a UFD and $D$ is a non-trivial irreducible $R$-LND without having a slice, by Corollary \ref{Cor_FibRk1-2LND}(I) and (II) we respectively have $\ResRk{D}=2$ and $\ResVarRk{D} =2$; and also by Theorem \ref{Rk2-LND_Properties_Pol-Alg} it follows that $\Rk{D} =2$. So, we have $\ResRk{D} = \ResVarRk{D} = \Rk{D} =2$. Note that since $A$ is an $\A{2}$-fibration over $R[F]$, by Lemma \ref{Lem_Fib-to-ResSys} we have $(R, R[F], A)$ is a $(3, 2)$-residual system such that $R[F] \subseteq \KerD$.
	\end{ex}
	
	One can notice that in each of examples \ref{Ex_Winkelmann-1} and \ref{Ex_BD_A2} the residual rank and residual-variable rank of the LNDs of $R^{[3]}$ are equal. We end this discussion by quoting an example by Raynaud (\cite{Raynaud}) as studied by Essen and Rossum (\cite{Essen-Rossum_counterexamlpes}) and Freudenburg (\cite{Freudenburg_Rxyz-slice}) (also see \cite{Freudenburg-BookNew}, Theorem 10.25) which shows that the residual rank and residual-variable rank of an LND of polynomial algebra need not agree, and also establish the fact that rank of an LND of a polynomial algebra $A$ need not be same as the rank of the trivial extension of the LND of $A^{[n]}$.
	
	\begin{ex} \label{Ex_Raynaud} \normalfont 
		Let $\displaystyle R = \mathbb{R}[a,b,c,x,y,z]/(ax + by + cz -1)$. Set $A := R[X,Y,Z]$. Define an LND $D:A \longrightarrow A$ by $D(X) = x$, $D(Y)= y$ and $D(Z) = z$. One can see that $D(aX + bY +cZ) =1$ and $\KerD = R [X -xs, \ Y-ys, \ Z-zs] \ne R^{[2]}$, and therefore, by Theorem \ref{Rk2-LND_Properties_Pol-Alg} we see that $\Rk{D} > 2$, i.e., $\Rk{D} =3$. Since $aX + bY + cZ$ is a slice of $D$, we have $\KerD^{[1]} = \KerD[aX + bY + cZ] = A = R^{[3]}$, and therefore, by Corollary \ref{Cor_FPF_LND_on_A^3} it follows that $\KerD$ is an $\A{2}$-fibration over $R$. Let $\tilde{D}: A[W] \longrightarrow A[W]$ be the trivial extension of $D$, i.e, $\tilde{D}(W) = 0$. Note that $\Ker{\tilde{D}}^{[1]} = A[W]$ and $\Ker{\tilde{D}} = \KerD[W] \cong \KerD[aX + bY + cZ] = \KerD^{[1]} = R^{[3]}$, and therefore, $\Rk{\tilde{D}} =1$.
		
		\smallskip
		
		Since $D$ has a slice, by Corollary \ref{Cor_FPF_LND_on_A^3} it follows that $\ResRk{D} =1$. We claim that $\ResVarRk{D} =3$. If possible, let $\ResVarRk{D} \le 2$. Then, by Corollary \ref{Cor_FibRk_FPF-LND}[(I) \& (III)] we have $\KerD = R^{[2]}$, which is a contradiction to the fact that $\KerD  \ne R^{[2]}$, and therefore, $\ResVarRk{D} =3$. So, we have $\ResRk{D} =1$ and $\ResVarRk{D} = \Rk{D} =3$.
	\end{ex}
	
	\section*{Acknowledgment:} The authors thank Amartya K. Dutta and Neena Gupta for helpful suggestions. 
	
	\bibliographystyle{apalike}
	\normalem
	\bibliography{reference}
\end{document}